\documentclass[12pt]{article}
\usepackage{amssymb,latexsym,amsmath}

\newtheorem{theorem}{Theorem}
\newtheorem{corollary}[theorem]{Corollary}
\newtheorem{lemma}[theorem]{Lemma}

\begin{document}

\title{On a Question of Arveson about Ranks of Hilbert modules}
 
\author{Xiang Fang}

\maketitle

%%%%%% Beginning of self-defined notations ! 

% 1 :  numbers

\def\cc{\mathbb{C}}
\def\ccc#1{\mathbb{C}^{#1}}

\def\zz{\mathbb{Z}}
\def\zzz#1{\mathbb{Z}^{#1}}

\def\nn{\mathbb{N}}
\def\nnn#1{\mathbb{N}^{#1}}

\def\rr{\mathbb{R}}
\def\rrr#1{\mathbb{R}^{#1}}

% 2 : disc algebras, polynomial ring
\def\ad#1{A({\mathcal{D}}^{#1})}
\def\ab#1{A({B}_{#1})}
\def\polyz{C[z]}
\def\polyzw{C[z,\omega]}
\def\poly#1{C[z_{1}, \cdots, z_{#1}]}

% 3 : spaces

\def\Hi{\mathcal{H}}
\def\K{\mathcal{K}}
\def\L{\mathcal{L}}

\def\htwo{\mathcal{H}^{2}} 
\def\htwoo#1{\mathcal{H}^{2} \otimes {\cc}^{#1} } 

\def\submo{\mathcal{M}}
\def\submoperp{{\submo}^{\perp}}

% 4 : operator tuple, coordinate tuple

\def\ztuple#1{{#1}_{z}=({#1}_{z_{1}}, \cdots, {#1}_{z_{d}}) }
\def\tuple#1{{#1}=({#1}_{1}, \cdots, {#1}_{d})}
\def\tup#1#2{{#1}_{1}, \cdots, {#1}_{#2}}

% 5 : Remark Proof and Section :

\def\remark{\noindent \emph{Remark}}
\def\proof{\noindent \emph{Proof}}
\def\sec#1{\noindent {\bf{#1}}} 

% 6 : range < ...< 

\def\range#1#2{1 \leqslant {#1} \leqslant {#2}}
\def\ran#1{1 \leqslant {#1} \leqslant d}

% 7 : sum 
\def\summ#1#2#3{\sum_{{#1}={#2}}^{{#3}}}
\def\sumd#1{\sum_{{#1}=1}^{d}}
\def\sumninfty{\sum_{n = 1}^{\infty}}
\def\sumna{\sum_{n = 1}^{\infty}{\alpha}_{n}^{-2}}

% 8 : Various

\def\en{\Box}

\def\tpr{T \cong Pr_{\submoperp}{(M_{z})}}
\def\tprN{T \cong Pr_{\submoperp}{(M_{z} \otimes I_{N})}}

\def\fred#1{Fredholm \quad index ({#1})}
\def\norm#1{\Vert{#1}\Vert}
\def\projection#1#2{P_{(#1)}(#2)}
\def\compression#1#2{Pr_{(#1)}({#2})}

%%%% This is the end of definitions ! %%%%

\begin{abstract}
It's well known that the functional Hilbert space $\htwo$ over the unit ball 
$B_{d}  \subset \ccc{d}$, with kernel function
$K(z,\omega)=\frac{1}{1-z_{1}\omega_{1}-\cdots -z_{d}\omega_{d}}$,  
admits a natural $\ab{d}$-module structure.
We show the rank of a nonzero submodule 
$\submo\subset\htwo$ is infinity if and only if $\submo$ is of infinite
codimension. Together with Arveson's dilation theory,
 our result
shows that Hilbert modules stand in stark contrast with 
 Hilbert basis theorem for
 algebraic modules.
This result answers a question of Arveson.
\end{abstract}

\bigskip

\noindent{$\quad\quad$\scriptsize{KEYWORDS: Hilbert module, submodule, rank.}}

\bigskip

\noindent {$\quad\quad$\scriptsize{MATHEMATICS SUBJECT CLASSIFICATION: 
           Primary 47A13; Secondary 47A20 46H25.}}

\bigskip

\sec{0. Introduction :} In the study of Multivariable Operator Thoery
$(MOT)$, there is a natural approach via Hilbert modules. 
(Readers are refered to \cite{DP} for more information.)
Let $\tuple{T}$ be a tuple of commuting operators acting on a common 
separable
Hilbert space $\Hi$. There is 
a natural $\ab{d}$-module structure on $\Hi$:
\begin{displaymath}
 f\cdot\xi = f(T_{1},\cdots, T_{d})\xi, \qquad f \in \ab{d},\quad \xi \in \Hi.
\end{displaymath}
This approach allows us to introduce algebraic 
techniques to operator theory. This 
is not so significant in one variable, 
since function theory in one variable
is often mature enough for analyzing problems. 
But in higher dimensions, the situation
is very different. We have to deal with algebraic and analytic
 varieties, which are discrete points in one
variable. Hence some algebraic tools, like localization, are useful and
indispensible. 
( See  \cite{berg}, \cite{DY}, \cite{DPSY}, \cite{Guo},$\cdots$ 
 for some nontrivial applications of 
algebraic techniques in operator
theory.)

The purpose of this note is to show a phenomenon of Hilbert modules 
which is in stark contrast with Hilbert basis theorem for finitely generated
modules over Neotherian rings.
 It suggests topological modules may behave very
differently from algebraic modules.

However in \cite{polynomial},  we will show that another aspect of algebraic 
 modules, namely the Hilbert polynomial, suits well 
 for Hilbert modules.

\bigskip

\sec{ 1. Preliminaries: a brief review of $\mathbf{\htwo}$ : } In this section, 
we review briefly the basic
facts about the Hilbert modules we are intersted in. Readers can find more
information in \cite{subalgebra}, \cite{summary}, \cite{curvature}.

Let $B_{d}  \subset \ccc{d}$ be the open unit ball in $\ccc{d}$.
Let $\htwo$ be the functional Hilbert space over $B_{d}  \subset \ccc{d} $,
determined by the reproducing kernel
$K(z,\omega)=\frac{1}{1-z_{1}\omega_{1}-\cdots - z_{d}\omega_{d}}$, where
$\tuple{z},\tuple{\omega} \in \ccc{d}$. In particular, different monomials 
in $\htwo$ are
orthogonal: $\langle z^{I}, z^{J} \rangle = 0, \quad I \ne J$, where 
 $I,J \in {\zz}_{d}^{+}$ are multi-indices.
Let $\ztuple{M}$ be the tuple of multiplication by coordinate functions on $\htwo$.
Then $M_{z}$ is a tuple of commuting, bounded operators. We have 
\begin{displaymath}
 I-M_{z_{1}}{M_{z_{1}}}^{*}-\cdots-M_{z_{d}}{M_{z_{d}}}^{*}=P_{0},
\end{displaymath}
where $P_{0}$ is the projection onto the constant term in $\htwo$.
Except for $P_{0}$, we will use $P_{(L)}$ denoting the 
orthogonal projection with the range  $L$.

For any tuple $\tuple{T}$ of commuting bounded operators on $\Hi$, following
Arveson's language, we say $T$ is $d$-contractive if 
\begin{displaymath}
I-T_{1}T_{1}^{*}-\cdots -T_{d}T_{d}^{*} \geq 0.
\end{displaymath}
So it follows $\ztuple{M}$ is $d$-contractive. For any 
$d$-contractive operator tuple $\tuple{T}$ on $\Hi$, the {\emph{``defect operator"}}
 is defined to be 
\begin{displaymath}
 \Delta=(I-T_{1}T_{1}^{*}-\cdots -T_{d}T_{d}^{*})^{\frac{1}{2}}.
\end{displaymath}
The rank of $\tuple{T}$ is defined to the dimension of $\overline{\Delta\Hi}$.
See \cite{subalgebra} for more information on this definition. By the dilation 
theory in \cite{subalgebra}, for any $d$-contractive tuple $\tuple{T}$ on $\Hi$,
 there exists an auxiliary Hilbert space $\L$, a spherical isometry $\tuple{S}$ on 
 $\L$, $(i.e.\quad S_{1}S_{1}^{*}+\cdots+S_{d}S_{d}^{*}=I)$, 
 and coinvariant subspaces $\K \subset \htwoo{(rank(T))}$ and 
 $\K' \subset \L$ such that
\begin{displaymath}
T \cong \compression{\K}{M_{z}\otimes I_{rank(T)}}\oplus \compression{\K'}{S},
\end{displaymath}
where the coinvariant subspaces are joint coinvariant with respect to obvious
actions. Since it is usually clear from the context which operators we
 are talking about, we will use the term ``$($co-$)$invariant" 
 without specifying the actions.
 When the second term $\compression{\K'}{S}$ is null, $T$ is called a
 \emph{``pure"} $d$-contraction.
 
Let $R_{z_{j}}=M_{z_{j}}\vert_{\submo}$ $(\ran{j})$ be the restrictions of 
$M_{z_{j}}$ $(\ran{j})$ to the submodule $\submo$. Because
\begin{eqnarray*}
\sumd{j} P_{\submo}M_{z_{j}}P_{\submo}\cdot
P_{\submo}M_{z_{j}}^{*}P_{\submo} & \leqslant &
\sumd{j} P_{\submo}M_{z_{j}}M_{z_{j}}^{*}P_{\submo} \\
\leqslant P_{\submo}P_{\submo} & = & P_{\submo},
\end{eqnarray*}
$\ztuple{R}$ is $d$-contractive. If the spherical isometry part of $R$ is
not null, it implies $R_{z}$ has a nontrivial reducing subspace. This reducing
subspace will also be a reducing subspace of $\ztuple{M}$ on $\htwo$. It is
impossible. So $\ztuple{R}$ is ``pure."

The rank of submodule $\submo$ is defined to be the rank of the tuple of
 restrictions $\ztuple{R}$. The main result of this note is 
 
\begin{theorem}[Main Result]\label{T:main}
    Let $\submo \subset \htwo$ be a nonzero submodule. Then $\submo$ has rank infinity
     if $\submo$ is of infinite codimension.
\end{theorem}

\remark :$(1)$ If  submodule $\submo$ is of finite codimension,
 it is easy to see
the rank of $\submo$ is finite. \\
$(2)$When $\submo$ is finitely generated by homogeneous 
polynomials, the above result is also proved in \cite{curvature}. \\
$(3)$ The above result should be compared
 with the corresponding result in commutative algebra.
Assume $R$ is a Noetherian ring. $($ $i.e.$ $R$ is 
a commutative ring with unit $e$ and any ideal of $R$ admits a finite 
system of generators. $)$ Then $F=R$ will be the free module of rank one.
Assume $M \subset F^{(n')}$ is a submodule of the free module
$F^{(n')}$. Then by Hilbert basis theorem, $M$ admits a 
 finite system of generators 
$\{g_{1},\cdots,g_{n}\}$. For the free module
 $F^{(n)}$ with finite rank $n$, there exists a unique surjective module 
 map $\varphi: F^{(n)} \to M$ which maps the canonical basis vector $e_{i}$ 
 to $g_{i}$ $(\range{i}{n})$. Let $N=ker(\varphi) \subset F^{(n)}$.  
Then $N$ is a submodule of the free module $F^{(n)}$ such that 
\begin{displaymath}
0 \to N \to F^{(n)} \to M \to 0
\end{displaymath}
is a short exact sequence. 

Recall that in \cite{curvature}, $\htwo$ is justified to be the free
module of rank one in the category of pure $d$-contractive Hilbert 
modules. By the dilation theory
in \cite{subalgebra}, if $\tilde{M}$ is a pure $d$-contractive Hilbert 
module with 
rank $\tilde{n}$, then there exists a submodule 
$\tilde{N}\subset \htwoo{(\tilde{n})}$, such that 
${\tilde{N}}^{\perp}$ is unitarily equivalent to $\tilde{M}$ as pure
 $d$-contractive modules. Or, we have the following short exact sequence
\begin{displaymath}
0 \to \tilde{N} \to \htwoo{(\tilde{n})} \to \tilde{M} \to 0.
\end{displaymath}

If $\tilde{M}\subset \htwo$ is a 
submodule of the free module $\htwo$ with infinite codimension, 
 by dilation theory in \cite{subalgebra} and the above theorem, 
 we cannot find
a finite integer $\tilde{n}$, such that $\tilde{M}$ is unitarily equivalent
to ${\tilde{N}}^{\perp}$ for some 
submodule $\tilde{N}\subset \htwoo{(\tilde{n})}$. Or, we cannot make the
above short sequence of Hilbert modules exact with some finite $\tilde{n}$.

\begin{corollary}\label{C:coro}
Let $\submo \subset \htwo$ be a nonzero submodule. Then
 there exists some integer $n \in \nn$ and some submodule 
  $\mathcal{N} \subset \htwoo{(n)}$, such that $\submo$ is  unitarily 
  equivalent to ${\mathcal{N}}^{\perp}$ as pure $d$-contractive Hilbert 
  modules if and only if $\submo$ is of finite codimension in $\htwo$.
\end{corollary}

\bigskip

\sec{2. Jumping operators $\mathbf{J_{z_{j}}:\submoperp \to \submo \quad (\ran{j})}$ :} 
 First we decompose $M_{z_{j}}$ $(\ran{j})$ with
  respect to $\submo$ and $\submoperp$ as :
\begin{displaymath}
M_{z_{j}} =
\left( \begin{array}{cc}
R_{z_{j}} & J_{z_{j}} \\
0 & S_{z_{j}} 
\end{array} \right), \quad (\ran{j}).
\end{displaymath}
Most studies of submodules and quotient modules center around $R_{z_{j}}$ and
$S_{z_{j}}$ $(\ran{j})$ so far. 
Less attention has been devoted to 
the ``Jumping operators" $J_{z_{j}}: \submoperp \to \submo$ $(\ran{j})$.
Our study in this note is closely connected 
to $J_{z_{j}}: \submoperp \to \submo$ $(\ran{j})$. In fact we find
 these operators have more significance when we
study Fredholm indices, or cohomology groups of a Koszul complex.

Observe that
\begin{eqnarray*}
P_{0} & = &  I - \sumd{j} M_{z_{j}}M_{z_{j}}^{*} \\
& = & \left( \begin{array}{cc} 
I - \sumd{j} R_{z_{j}}R_{z_{j}}^{*} - \sumd{j} J_{z_{j}}J_{z_{j}}^{*} & \quad\cdots \\
\cdots & \quad\cdots \end{array} \right).
\end{eqnarray*}
We use $\compression{\submo}{P_{0}}$ denoting the compression of $P_{0}$ 
onto subspace $\submo$.
Since $\compression{\submo}{P_{0}}$ is of rank one, 
we know $rank(I-\sumd{j} R_{z_{j}}R_{z_{j}}^{*}) < \infty$ if and only if 
$rank(\sumd{j} J_{z_{j}}J_{z_{j}}^{*}) < \infty$. This is equivalent to say
$rank(J_{z_{j}}) < \infty$ for every $j$ $(\ran{j})$. So it suffices
 to prove the following
characterization of infinite codimensional submodule
 $\submo \subset \htwo$.

\begin{theorem}\label{T:theorem2}
Let $\submo \subset \htwo$ be a nonzero submodule. Then 
$\submo$ is of infinite codimension  if and only if 
\begin{displaymath}
rank(J_{z_{1}})+\cdots + rank(J_{z_{d}}) = \infty
\end{displaymath}
\end{theorem}

\remark : We will also prove the above result for the Hardy
space $\htwo(\mathbb{D}^{d})$ in the last section of this note. The proof relies on
Beurling's theorem and the fact that $\ztuple{R}$ are all isometries.
 But it seems the idea for the Hardy space
$\htwo(\mathbb{D}^{d})$ will carry over to $\htwo$, provided 
we have a good understanding of the invariant subspaces of $\htwo$. It is
interesting to compare the proof of analogous results for 
$\htwo(\mathbb{D}^{d})$ and $\htwo$. Since our ultimate goal is to develop a
theory for commutting operator tuples, we would 
like to see a proof of the main result in 
this note via
analyzing the behavior of invariant subspaces of $\htwo$.
Especially, we would like to see some results for $\htwo$, which are
 analogous 
 to Lax-Halmos theorem for Hardy spaces.

\bigskip

\sec{3 : Two lemmas :} We introduce two lemmas needed 
for the next section. Moreover, we will use these two lemmas for 
the Hardy space $\htwo(\mathbb{D}^{d})$ in the last section without further 
explanation since the proofs here carry over to 
the Hardy space $\htwo(\mathbb{D}^{d})$.

\begin{lemma}\label{L:lemma1}
 $dim(\submo \ominus z_{j}\submo) = \infty$, \quad   ($\ran{j})$.
\end{lemma}

We need to introduce some notation before we give the proof. 
For any function $f \in \htwo$, let
$f=\summ{n}{0}{\infty} f^{(n)}$ be the homogeneous expansion of $f$, such that 
$f^{(n)}$ is a homogeneous polynomial of degree $n$. We define
\begin{displaymath}
  ord(f)=inf\{n: \quad f^{(n)}\ne 0,\quad f=\summ{n}{0}{\infty} f^{(n)} \}.
\end{displaymath}
For a submodule $\submo \subset \htwo$, we  define
\begin{displaymath}
  ord(\submo)=inf\{ord(f): \quad f \in \submo\}.
\end{displaymath}
Let $I=(\tup{z}{d})$ be the maximal idea of $\ad{d}$ at the origin. Then 
it is easy to see $ord(\overline{I^{k}\submo})=ord(\submo)+k$ $(k \in \nn)$ 
$($ since $ord(\cdot)$ is upper semi-continous$)$.

\proof : We fix some $j$ $(\ran{j})$ first.

 Assume $dim(\submo\ominus z_{j}\submo) = k < \infty$. We first show
\begin{displaymath}
dim(\submo\ominus z_{j}^{l}\submo) \leqslant l\cdot k, \qquad l \in \nn.
\end{displaymath}
 It is easy to see
\begin{eqnarray*}
dim(\submo\ominus z_{j}^{l}\submo) & =  & dim(\submo\ominus z_{j}\submo) + 
dim(z_{j}\submo\ominus z_{j}^{2}\submo) + \\
& & \qquad \qquad\quad\cdots+dim(z_{j}^{l-1}\submo\ominus z_{j}^{l}\submo).
\end{eqnarray*}
For each $t=1,2,\cdots,l$, we define 
$B_{t}:z_{j}^{t-1}\submo\ominus z_{j}^{t}\submo \to \submo\ominus z_{j}\submo$
by
\begin{displaymath}
B_{t}(z_{j}^{t-1}x)=\projection{\submo\ominus z_{j}\submo}{x},\qquad 
z_{j}^{t-1}x \in z_{j}^{t-1}\submo\ominus z_{j}^{t}\submo,
\end{displaymath}
where $x \in \submo$.

 If $B_{t}(z_{j}^{t-1}x) = 0$, $x \in \submo$ and 
$\projection{\submo\ominus z_{j}\submo}{x}=0$ will imply $x \in z_{j}\submo$.
Hence $ z_{j}^{t-1}x \in z_{j}^{t}\submo$. It follows 
that $B_{t}$ is injective,
and $dim(\submo\ominus z_{j}^{l}\submo) \leqslant lk$.

Assume $ord(\submo)=a$. Then $ord(f)=a$ for some $f\in \submo$. We claim
\begin{displaymath}
  dim(\submo\ominus I^{l}\submo) \geqslant C_{l+d-1}^{d}, \qquad l\in \nn.
\end{displaymath}
$($In the above we may also write 
$\submo\ominus\overline{I^{l}\submo}$ since $I^{l}\submo$ is not generally 
closed. But we want to make the notation simple and do not do that.$)$ 
For a fixed $l \in \nn$, we define 
\begin{displaymath}
L=span\{z^{I}\cdot f:\quad |I|=i_{1}+\cdots+i_{d}<l\}.
\end{displaymath}
Then $P_{(\submo\ominus I^{l}\submo)}\vert_{L}:L\to\submo\ominus I^{l}\submo$
 is injective.
 
  Otherwise, assume
$\projection{\submo\ominus I^{l}\submo}{p\cdot f}=0$ for some polynomial $p$
 with $deg(p) < l$. Then $ord(p\cdot f)\leqslant l+a-1$. But   
$\projection{\submo\ominus I^{l}\submo}{p\cdot f}=0$ implies
$p\cdot f \in \overline{I^{l}\submo}$. 
So $ord(p\cdot f)\geqslant l+a$. It is impossible.

So it follows 
\begin{displaymath}
dim(\submo\ominus I^{l}\submo) \geqslant dim(L)
\geqslant C_{l+d-1}^{d}, \quad l\in \nn.
\end{displaymath}

Because $z_{j}^{l}\submo \subset I^{l}\submo$, we conclude  
\begin{displaymath}
C_{l+d-1}^{d}\leqslant dim(\submo\ominus I^{l}\submo) \leqslant 
dim(\submo\ominus z_{j}^{l}\submo) \leqslant l\cdot k, \qquad l \in \nn.
\end{displaymath}
But when $l$ is large enough, it forces $k= \infty$. $\en$

\bigskip

\remark : By a remarkable theorem due to David Hilbert and \cite{DY}, the
 dimension function 
$\varphi(l)=dim(\submo\ominus I^{l}\submo)$ will become a polynomial 
when $l$ is large enough.
 $($See \cite{DY}, \cite{hart} for more information.$)$ 
This polynomial, named after Hilbert, is well-known 
to algebraists, but less familar to the operator theorists. 
In the Hilbert module setting, it has some nice properties and applications.
 We will pursue this topic in \cite{polynomial}.

Let $\Hi_{z_{i{1}},\cdots,z_{i{k}}} \subset \htwo$ 
$(1 \leqslant z_{i{1}} < \cdots < z_{i{k}} \leqslant d)$ be the subspace of
$\htwo$, consisting of elements involving only constants and variables
$z_{i{1}},\cdots,z_{i{k}}$. Then $\Hi_{z_{1}},\cdots,\Hi_{z_{d}}$ are the
classical Hardy spaces with variables $z_{1},\cdots,z_{d}$ respectively.

\begin{lemma}\label{L:lemma2}
A submodule $\submo \subset \htwo$ is of finite codimension in $\htwo$ if and
only if  $\submo \cap \Hi_{z_{j}} = u_{j}(z_{j})\Hi_{z_{j}}$ for some finite
nonzero  Blaschke product $u_{j}(z_{j})$ in variable $z_{j}$, $(\ran{j})$.
\end{lemma}

\proof : If $dim(\submoperp)< \infty$, the codimension of 
$\submo \cap \Hi_{z_{j}}$ in $\Hi_{z_{j}}$ is finite for each $j$
$(\ran{j})$. It is well known that
an invariant subspace of the classical Hardy space is of
 finite codimension if
and only if the corresponding inner function is a finite Blaschke product. 

If $\submo \cap \Hi_{z_{j}} = u_{j}(z_{j})\Hi_{z_{j}}$ for some finite
 Blaschke product $u_{j}$ in variable $z_{j}$ for each $j$
$(\ran{j})$, 
$\submo \cap \Hi_{z_{j}}$ then contains a nonzero polynomial $p_{j}(z_{j})$.
Assume $deg(p_{j}(z_{j})) \leqslant A$ for each $j$ $(\ran{j})$. Let
$\pi : \htwo \to \htwo/\submo = \submoperp$ be the quotient map. Then 
\begin{displaymath}
\submoperp=span\{\pi(z^{I}) : 
\quad I=(i_{1},\cdots, i_{d}),i_{r} < A, 1\leqslant r \leqslant d\}
\end{displaymath} 
is of finite dimension. $\en$

\bigskip

\sec{ 4. Proof of main result for ${\mathbf{d=2}}$ : }
We will prove the main result for $d=2$, and then use 
induction argument to finish the general case.
As we have mentioned, since our goal is to develop $MOT$, we 
would like to see a proof without using induction argument.

Let $\K=\submoperp$ and assume $dim(\K)=\infty$. 
For each $j$ $(\ran{j})$, we have 
\begin{eqnarray*}
{\{\K+z_{j}\K+z_{j}\submo\}}^{\perp} & = & 
             \submo \cap \Hi_{z_{1},\cdots,\hat{z_{j}},\cdots,z_{d}} \\
     & = & (\submo\ominus z_{j}\submo) \cap \Hi_{z_{1},\cdots,\hat{z_{j}},\cdots,z_{d}}.
\end{eqnarray*}
We decompose  $\submo\ominus z_{j}\submo$ $(\ran{j})$ into orthogonal sums as 
\begin{displaymath}
\submo\ominus z_{j}\submo=
[(\submo\ominus z_{j}\submo) \cap \Hi_{z_{1},\cdots,\hat{z_{j}},\cdots,z_{d}}]
\oplus \mathcal{E}_{j},
\end{displaymath}
where $\mathcal{E}_{j} \subset \overline{\K+z_{j}\K+z_{j}\submo}$, $\ran{j}$.
 So we have 
\begin{displaymath}
\mathcal{E}_{j} \subset 
   \projection{\submo\ominus z_{j}\submo}{z_{j}\K},\quad \ran{j}.
\end{displaymath}   
For any $x\in \mathcal{E}_{j}$ $(\ran{j})$, we may write 
$x=\projection{\submo\ominus z_{j}\submo}{z_{j}\xi}$ $(\ran{j})$,
 where $\xi \in \K$.
Then $\projection{\submo\ominus z_{j}\submo}{J_{z_{j}}(\xi)}=x$
 $(\ran{j})$. It follows 
\begin{displaymath}
rank(J_{z_{j}}) \geqslant dim(\mathcal{E}_{j}), \qquad \ran{j}.
\end{displaymath}

From now on, we assume 
\begin{displaymath}
dim(\mathcal{E}_{j}) < \infty,\qquad \ran{j}.
\end{displaymath}
This assumption will be used as we finish our proof by induction in the next section.

For the rest of this section, we assume $d=2$.

By Beurling's theorem, there is an inner function 
     $\varphi(z_{2}) \in \Hi_{z_{2}}$, such that
\begin{eqnarray*}
(\submo\ominus z_{1}\submo) \cap \Hi_{z_{2}}&=&\submo\cap \Hi_{z_{2}}\\
& = & \varphi(z_{2})\Hi_{z_{2}}.
\end{eqnarray*}
By Lemma~\ref{L:lemma1} $dim(\submo\ominus z_{1}\submo) = \infty$
 and our assumption $dim(\mathcal{E}_{1}) < \infty$, we know $\varphi(z_{2})$
 is nonzero.
 
For any polynomial $p(z_{1})$ in variable $z_{1}$, we define
\begin{displaymath}
S_{\{p\}}=\{g\in \submo: g=p(z_{1})f(z_{2}),  f(z_{2})\in Hol(\mathbb{D}) \}.
\end{displaymath}
Let $ord_{z_{2}}(f(z_{2}))$ be 
the multiplicity of the zeros at the origin of the 
holomorphic function $f(z_{2})$ in variable $z_{2}$.
For any $h(z_{1},z_{2})=g(z_{1})f(z_{2})$, where $g(z_{1}), f(z_{2})$
are holomorphic, we define
\begin{displaymath}
ord_{z_{2}}(h(z_{1},z_{2}))=ord_{z_{2}}(f(z_{2})).
\end{displaymath}
For polynomial $p(z_{1})$, we define
\begin{displaymath}
ord_{z_{2}}(S_{\{p\}})=inf\{ord_{z_{2}}(x) : \quad x\in S_{\{p\}}\},
\end{displaymath}
if $S_{\{p\}} \ne \{0\}$. For $S_{\{p\}} = \{0\}$, we define
$ord_{z_{2}}(S_{\{p\}})=+\infty$.
Let
\begin{displaymath}
a=inf\{ord_{z_{2}}(S_{\{p\}})  : \quad p \in C[z_{1}] \},
\end{displaymath}
where $C[z_{1}]$ is the polynomial ring in variable $z_{1}$. Since 
\begin{displaymath}
S_{\{1\}}=\submo\cap \Hi_{z_{2}} = \varphi(z_{2})\Hi_{z_{2}},
\end{displaymath}
we know $a<+\infty$.
Assume $p_{0}(z_{1})f_{0}(z_{2})$  achieves the value $a$.

Now we check that for any nonzero polynomial $p(z_{1})$ in $z_{1}$,
\begin{displaymath}
p(z_{1})p_{0}(z_{1})f_{0}(z_{2}) \notin z_{2}\submo.
\end{displaymath}
Otherwise, for some $p(z_{1})$, 
\begin{displaymath}
p(z_{1})p_{0}(z_{1})f_{0}(z_{2})=z_{2}x, \quad x \in \submo.
\end{displaymath}
So $f_{0}(z_{2})=z_{2}\tilde{f_{0}}(z_{2})$ for some holomorphic function
$\tilde{f_{0}}(z_{2})$ and 
$x=p(z_{1})p_{0}(z_{1})\tilde{f_{0}}(z_{2}) \in \submo$. Then 
$ord_{z_{2}}(x)=a-1$. Contradiction.

Recall that 
\begin{displaymath}
\submo\ominus z_{2}\submo=
[(\submo\ominus z_{2}\submo) \cap \Hi_{z_{1}}] \oplus \mathcal{E}_{2},
\end{displaymath}
where $dim(\mathcal{E}_{2})<\infty$.
Let 
\begin{displaymath}
\mathcal{G}=\{p(z_{1})p_{0}(z_{1})f_{0}(z_{2}):\quad  p(z_{1}) \in C[z_{1}] \}.
\end{displaymath}
Then $P_{(\mathcal{E}_{2})}\vert_{\mathcal{G}}:{\mathcal{G}}\to \mathcal{E}_{2}$ 
has finite dimensional range. Assume 
\begin{displaymath}
dim(P_{(\mathcal{E}_{2})}(\mathcal{G}))=b\leqslant dim(\mathcal{E}_{2}) <\infty.
\end{displaymath}
So there exists an $N\in \nn$ such that 
\begin{displaymath}
dim(P_{(\mathcal{E}_{2})}
(\{p(z_{1})p_{0}(z_{1})f_{0}(z_{2}):\quad 
 p(z_{1}) \in C[z_{1}], deg(p(z_{1}))\leqslant N  \}))=b.
\end{displaymath}
Hence for any $p(z_{1}) \in C[z_{1}]$, with $deg(p(z_{1})) > N$, there
exists some $\tilde{p}(z_{1}) \in  C[z_{1}]$, with  $deg(\tilde{p}(z_{1}))\leqslant N$,
such that
\begin{displaymath}
P_{(\mathcal{E}_{2})}((p-\tilde{p})p_{0}f_{0})=0.
\end{displaymath}
This means
\begin{displaymath}
(p-\tilde{p})p_{0}f_{0} \in [(\submo\ominus z_{2}\submo) \cap \Hi_{z_{1}}]\oplus z_{2}\submo.
\end{displaymath}
Hence there exists $\xi(z_{1}) \in (\submo\ominus z_{2}\submo) \cap \Hi_{z_{1}}$ and 
$\eta(z_{1},z_{2}) \in \submo$, such that
\begin{displaymath}
(p(z_{1})-\tilde{p}(z_{1}))p_{0}(z_{1})f_{0}(z_{2})=\xi(z_{1})+ z_{2}\eta(z_{1},z_{2}).
\end{displaymath}
Let $z_{2}=0$, we know  
\begin{displaymath}
\xi(z_{1}) = (p(z_{1})-\tilde{p}(z_{1}))p_{0}(z_{1})f_{0}(0)
\end{displaymath}
is a polynomial.
But we know
\begin{displaymath}
(p(z_{1})-\tilde{p}(z_{1}))p_{0}(z_{1})f_{0}(z_{2}) \notin z_{2}\submo.
\end{displaymath}
So $\xi(z_{1}) \in (\submo\ominus z_{2}\submo) \cap \Hi_{z_{1}}$ is a nonzero polynomial. 

Similarly, we can show $\submo$ contains a nonzero polynomial in $z_{2}$. Together with
Lemma~\ref{L:lemma2}, we finish our proof for $d=2$. $\en$

\bigskip

\sec{5. Proof of main result for {$\mathbf{d \geqslant 3}$ : }} Recall again
\begin{displaymath}
\submo\ominus z_{j}\submo=
[(\submo\ominus z_{j}\submo) \cap \Hi_{z_{1},\cdots,\hat{z_{j}},\cdots,z_{d}}]
\oplus \mathcal{E}_{j},
\end{displaymath}
where we assume $dim(\mathcal{E}_{j}) < \infty$ $(\ran{j})$.

Assume our proof is finished for $(d-1)$ dimensional case.

If for some $j$ $(\ran{j})$, 
$\submo \cap \Hi_{z_{1},\cdots,\hat{z_{j}},\cdots,z_{d}}$ is 
of finite codimension in $\Hi_{z_{1},\cdots,\hat{z_{j}},\cdots,z_{d}}$, then 
$\submo \cap \Hi_{z_{s}}$ is of finite codimension in $\Hi_{z_{s}}$ 
for any $s$ such that $\ran{s}$, $s\ne j$. 

So if for any $j$ $(\ran{j})$, 
$\submo \cap \Hi_{z_{1},\cdots,\hat{z_{j}},\cdots,z_{d}}$ is 
of finite codimension in $\Hi_{z_{1},\cdots,\hat{z_{j}},\cdots,z_{d}}$,
Lemma~\ref{L:lemma2} will imply $\submo$ is of finite codimension in $\htwo$.

Now we assume for some $j$ $(\ran{j})$, say $j=d$, 
$\submo \cap \Hi_{z_{1},\cdots,z_{d-1}}$ is 
of infinite codimension in $\Hi_{z_{1},\cdots,z_{d-1}}$. Assumption
$dim(\mathcal{E}_{j}) < \infty$ $(\ran{j})$ 
implies $\submo \cap \Hi_{z_{1},\cdots,z_{d-1}}$
is a nonempty submodule in $\Hi_{z_{1},\cdots,z_{d-1}}$.
Now we may apply induction.

Our proof of the main result is completed. $\en$

\bigskip

\sec{ 6. A analogous result for the Hardy space $\mathbf{\htwo(\mathbb{D}^{d})}$ : }
In this section, we will prove  Theorem~\ref{T:theorem2} for
Hardy space $\htwo(\mathbb{D}^{2})$
 over the bidisk $\mathbb{D}^{2}$. The reason why we choose
bidisk $\mathbb{D}^{2}$ instead of the polydisk $\mathbb{D}^{d}$ is only for 
notational convenience. The proof here relies heavily on Beurling's theorem and
the fact that $\ztuple{R}$ are isometries. Let $(z,\omega)$ be the coordinate
functions on  $\mathbb{D}^{2}$. $\Hi_{z}$ and $\Hi_{\omega}$ have the same
meaning as in Section $3$.
We decompose $M_{z_{j}}$ $(\ran{j})$ with respect to 
a submodule $\submo$ and $\submoperp$ in the same way as we did 
in Section $2$:
\begin{displaymath}
M_{z_{j}} =
\left( \begin{array}{cc}
R_{z_{j}} & J_{z_{j}} \\
0 & S_{z_{j}} 
\end{array} \right), (\ran{j}).
\end{displaymath}
We will prove 

\begin{theorem}\label{T:theorem3}
A nonzero submodule $\submo \subset \htwo(\mathbb{D}^{d})$
 is of infinite codimension  if and only if 
\begin{displaymath}
rank(J_{z_{1}})+\cdots + rank(J_{z_{d}}) = \infty
\end{displaymath}
\end{theorem}

\proof : 
In fact we will deal with $J_{z}^{*}$ and $J_{\omega}^{*}$, instead of
$J_{z}$ and $J_{\omega}$. First assume $\submo$ is of infinite codimension.

 Observe that $J_{z}^{*}(z\cdot x)=0$ for any 
$x \in \submo$, so we only have to consider the rank of 
\begin{displaymath}
J_{z}^{*} : \submo\ominus z\submo \to \submoperp.
\end{displaymath}
First we look at the kernel of 
$J_{z}^{*}$ on $\submo\ominus z\submo$.
Let 
\begin{displaymath}
L_{z}=\{x \in  \submo\ominus z\submo:\quad  J_{z}^{*}(x)=0 \} 
\end{displaymath}
Then for any $x \in L_{z}$, $i.e.$ $P_{\submoperp}(\overline{z}x)=0$, 
we have 
$P_{\htwo}(\overline{z}x) \in \submo$. But for any 
$x \in \submo\ominus z\submo$, we have
$\langle \overline{z}x, \submo \rangle =0$. It follows 
$P_{\htwo}(\overline{z}x)=0$, and $x \in \Hi_{\omega}$. It is easy to 
see that if $x \in \Hi_{\omega} \cap \submo$, then $x \in L_{z}$. So it follows
$L_{z}=\Hi_{\omega} \cap \submo$.

If $rank(J_{z}^{*}) < \infty$, then there exists a finite dimensional
subspace $\tilde{L_{z}} \subset \submo\ominus z\submo$ such that 
\begin{eqnarray*}
\submo\ominus z\submo & = & L_{z}+ \tilde{L_{z}} \\
& = & \varphi(\omega)\Hi_{\omega} + \tilde{L_{z}} .
\end{eqnarray*}
Similarly, if $rank(J_{z}^{*}) < \infty$, we have 
\begin{eqnarray*}
\submo\ominus \omega\submo & = & L_{\omega}+ \tilde{L_{\omega}} \\
& = & \psi(z)\Hi_{z} + \tilde{L_{\omega}}.
\end{eqnarray*}
Here $\varphi(\omega)$ and $\psi(z)$ are inner functions.
Because of Lemma~\ref{L:lemma2}, we may assume 
$\varphi(\omega)\Hi_{\omega}$ is of infinite codimension in 
$\Hi_{\omega}$. 

Consider $\submo$ is a $z$-invariant subspace of 
$\htwo(\mathbb{D}^{2})=\Hi_{z}\otimes \Hi_{\omega}=
\htwo(\mathbb{D},\Hi_{\omega})$.
By Lax-Halmos theorem, we can find a subspace 
$\mathcal{E} \subset \Hi_{\omega}$ with infinite codimension in 
$\Hi_{\omega}$, such that the $z$-invariant subspace
$[\mathcal{E}]_{z} \subset \htwo(\mathbb{D}^{2})$, 
generated by $\mathcal{E}$, contains
the $z$-invariant subspace 
$[\tilde{L_{z}}]_{z} \subset \htwo(\mathbb{D}^{2})$, generated by 
$\tilde{L_{z}}$, and contains $L_{z}=\varphi(\omega)\Hi_{\omega}$. 
 Choose any inner function
$u(\omega) \in \Hi_{\omega}\ominus\mathcal{E}$. Then 
$u(\omega)\Hi_{z}$ is orthogonal to $[\mathcal{E}]_{z}$,
 hence orthogonal to $\submo$ since $\submo\subset [\mathcal{E}]_{z}$.
But $u(\omega)\psi(z) \in \psi(z)\Hi_{\omega} \subset \submo$.
 Contradiction. 
 
  The other direction is trivial. $\en$

\bigskip

\emph{Acknowledgements:} This work was done while the author 
was a graduate student at TAMU
 under the guidance of Dr. R. Douglas. 
 The author would like to thank his
 advisor for many valuable conversations and support.

\noindent Author's address:\\
Xiang Fang\\
 Department of Mathematics\\
Texas A\&M University\\
College Station, TX 77843\\
 USA
 
\bigskip 
\noindent 
E-mail: xfang@math.tamu.edu


\begin{thebibliography}{99} 
\bibitem{subalgebra}
   Arveson, W. :
   \emph{Subalgebras of ${C}^{*}$-algebras III:
           Multivariable operator theory,}
   Acta Math
   \textbf{181} (1998), 159-228.
\bibitem{summary}
   Arveson, W. :
   \emph{The Curvature of a Hilbert module over
            $\poly{d}$,}
   Proc. Nat. Acad. Sci. (USA) 
   \textbf{96} (1999),  11096-11099.
\bibitem{curvature}
   Arveson, W. :
   \emph{The Curvature invariant of a Hilbert
    module over $\poly{d}$,} 
    J. Reine Angew. Math. 
    \textbf{522} (2000), 173--236.  
\bibitem{DP}
  Douglas, R. and Paulsen, V. :
  \emph{Hilbert Modules over Function Algebras,}
  Pitman Research Notes in Mathematics
  \textbf{217}
  Longman Scientific \& Technical
  publaddr Harlow, Essex, UK.   
\bibitem{berg}
   Douglas, R. and Yan, K. :
   \emph{A multi-variable Berger-Shaw theorem,}
   J. Operator Theory 
   \textbf{27} (1992), no.1, 205--217.       
\bibitem{DY}
  Douglas, R. and Yan, K. :
  \emph{Hilbert-Samuel polynomials for Hilbert modules,}
  Indiana Univ. Math. J.
  \textbf{42},  (1993), no. 3, 811--820.
\bibitem{DPSY}
   Douglas, R., Paulsen, V., Sah, C.-H. and Yan, K. :
   \emph{Algebraic reduction and rigidity for Hilbert modules,}
   Amer. J. Math.
   \textbf{117} (1995), no. 1, 75--92.
\bibitem{polynomial}
  Fang, X. :
  \emph{A note on the 
  Hilbert polynomials of Hilbert modules,}  
  in preparation.
\bibitem{Guo}
   Guo, K. :
   \emph{Equivalence of Hardy submodules generated by polynomials,}
   J. Funct. Anal. 
   \textbf{178} (2000), no.2, 343--371.
\bibitem{hart}
   Hartshorne, R. :
   \emph{Algebraic geometry,}
   Graduate Texts in Mathematics, No. \textbf{52},
    Springer-Verlag, New York-Heidelberg, 1977.
\end{thebibliography}
\end{document}